\documentclass[a4paper,12pt,leqno]{article}

\usepackage{amssymb,latexsym}
\usepackage{mltex}
\usepackage{amsmath}
\usepackage{fancyhdr}
\usepackage{enumerate}
\usepackage{amsthm}
\newtheorem{thm}{Theorem}     
     
\newtheorem{lem}{Lemma}[section]
\newtheorem{prop}[lem]{Proposition}     
\newtheorem{rem}{Remark}     
\newtheorem{defn}[lem]{Definition}     
     
\newtheorem*{thrm}{Theorem}
\newcommand{\Ekt}{\mathbb{E}(\kappa,\tau)}
\newcommand{\Rbar}{\overline{R}}
\newcommand{\qphi}{Q_{\varphi}}

\newcommand{\nablab}{\overline{\nabla}}

\newcommand{\diff}{\mathrm{d}\,}

\newcommand{\iid}{\mathrm{Id}\,}

\newcommand{\ddet}{\mathrm{det}\,}
\newcommand{\ddiv}{\mathrm{div}\,}

\newcommand{\sspec}{\mathrm{spec}}

\newcommand{\trace}{\mathrm{tr\,}}

\newcommand{\lto}{\ensuremath{\longrightarrow}}
\newcommand{\C}{\mathbb{C}}
\newcommand{\Dt}{\dfrac{\partial }{\partial t}}
\newcommand{\dt}{\frac{\partial }{\partial t}}
\newcommand{\Ss}{\mathbb{S}}
\newcommand{\HH}{\mathbb{H}}

\newcommand{\R}{\mathbb{R}}

\newcommand{\M}{\mathbb{M}}

\newcommand{\Cinf}{\mathcal{C}^{\infty}}
\newcommand{\Chi}{\mathfrak{X}}
\newcommand{\pre}{\Re e}         

\newcommand{\base}{\{e_1,\ldots ,e_n\}}

\newcommand{\MR}{\mathbb{M}^2(\kappa)\times\R}

\newcommand{\MRn}{\mathbb{M}^n\times\R}
\newcommand{\function}[5]
{\begin{eqnarray*}\begin{array}{r@{}ccl}
 #1\;\colon\;  & #2 &\lto & #3 \\[.05cm]  
  & #4 &\longmapsto  & #5 
\end{array}\end{eqnarray*}
}
\newcommand{\beqt}{\begin{equation}}  \newcommand{\eeqt}{\end{equation}}
\newcommand{\bal}{\begin{align}}      \newcommand{\eal}{\end{align}}
\newcommand{\ba}{\begin{array}}      \newcommand{\ea}{\end{array}}
\newcommand{\bc}{\begin{center}}     \newcommand{\ec}{\end{center}}
\newcommand{\be}{\begin{enumerate}}  \newcommand{\ee}{\end{enumerate}}
\newcommand{\beq}{\begin{eqnarray}}  \newcommand{\eeq}{\end{eqnarray}}
\newcommand{\beQ}{\begin{eqnarray*}} \newcommand{\eeQ}{\end{eqnarray*}}
\newcommand{\bi}{\begin{itemize}}    \newcommand{\ei}{\end{itemize}}
\newcommand{\bt}{\begin{tabular}}    \newcommand{\et}{\end{tabular}}
\newcommand{\finpreuve}{\hfill\square\\}
\newcommand{\pfend}{\hfill\square}
\def\pf{\noindent{\textit {Proof :} }}

\title{Spinorial characterizations of Surfaces into 3-dimensional homogeneous Manifolds}
\author{Julien Roth}
\date{\today}
\begin{document}

\maketitle
\begin{center}
Institut \'Elie Cartan, UMR 7502\\
Nancy-Universit\'e, CNRS, INRIA\\
B.P. 239, 54506 Vand\oe uvre l\`es Nancy Cedex, France
\end{center}
\begin{center}
roth@iecn.u-nancy.fr
\end{center}
\begin{abstract}
We give a spinorial characterization of isometrically immersed surfaces into $3$-dimensional homogeneous manifolds with $4$-dimensional isometry group in terms of the existence of a particular spinor, called generalized Killing spinor. This generalizes results by T. Friedrich \cite{Fr} for $\R^3$ and B. Morel \cite{Mo} for $\Ss^3$ and $\HH^3$. The main argument is the interpretation of the energy-momentum tensor of a genralized Killing spinor as the second fondamental form up to a tensor depending on the structure of the ambient space.
\end{abstract}
{\it keywords:} Dirac Operator, Killing Spinors, Isometric Immersions, Gauss and Codazzi Equations.\\\\
\noindent
{\it subjclass:} Differential Geometry, Global Analysis, 53C27, 53C40, 
53C80, 58C40.

\date{}
\maketitle\pagenumbering{arabic}
\section{Introduction}
\label{sect1}

Over the past years, the spinorial tool has been used succesfully in the study of the geometry and the topology of submanifolds of space forms. The spinorial approach allows to solve naturally some problems of geometry of submanifolds. For instance, some simple proofs of the Alexandrov theorem in the Euclidean space were given (\cite{HMZ1}) or in the hyperbolic space (\cite{HMR}) and new results about Einstein manifolds or Langrangian submanifolds of K\"ahler manifolds (\cite{HMU}). \\
\indent
In the same time, several results about the $3$-dimensional homogeneous manifolds with $4$-dimensional isometry group were obtained, particularly concerning the minimal and constant mean curvature surfaces in this spaces, by H. Rosenberg {\it and al.} (\cite{AR,Ro} for instance). In this paper, we make use of the spinorial tool in the study of surfaces in these $3$-dimensional homogeneous manifolds. As in the case of space forms $\R^n$, $\Ss^n$ and $\HH^n$, we can think that a spinorial approach could yield to solve some open questions in the theory of surfaces in $3$-dimensional homogeneous manifolds with $4$-dimensional isometry group, as the existence of an Alexandrov type theorem.\\
\indent
The first step of such an approach is to understand the spinorial geometry of these surfaces and particularly to give a spinorial characterization of the surfaces which are isometrically immersed into these homogeneous manifolds. Precisely, we characterize these surfaces by the existence of a special spinor field called {\it generalized Killing spinor}. Then, we will show that the existence of a generalized Killing spinor is equivalent to the existence of a spinor field solution of a weaker equation (involving the Dirac operator) with an additional condition on the norm of the spinor field (see Theorems \ref{thm1}, \ref{thm2} and \ref{thm3}).\\
\indent
The results that we give in this paper are non-trivial generalizations of (\cite{Fr2}) and (\cite{Mo}) to these homogeneous manifolds.
\section{Preliminaries}
\label{sect2}
We begin with a preliminary section in order to recall the basic facts about spin geometry of hypersurfaces. For further details, the reader can refer to \cite{BHMM}, \cite{Fr3} or \cite{LM} for general properties of spin manifolds and \cite{Ba2}, \cite{HMo} or \cite{Mo3} for the restriction to hypersurfaces. 

\subsection{Hypersurfaces and induced spin structures}
\label{sect22}
Let $(M^{n+1},g)$ be a Riemannian spin manifold and $\Sigma M$ its spinor bundle. We denote by $\nabla$ the Levi-Civita connection on $TM$, and also on $\Sigma M$. The Clifford multiplication will be denoted by $\gamma$ and $\left\langle .,.\right\rangle $ is the natural Hermitian product $\Sigma M$, compatible with $\nabla$ and $\gamma$. Finally, we denote by $D$ the Dirac operator.\\
\indent
Now let $N$ be an orientable hypersurface of $M$. Since the normal bundle is trivial, the hypersurface $N$ is also spin. Indeed, the existence of a normal unit vector field $\nu$ globally defined induces on $N$ a spin structure from the spin structure of $M$.\\
\indent
Then we can consider the intrinsic spinor bundle $\Sigma N$ of $N$. We denote  respectively by $\nabla^{N}$,
$\gamma^{N}$ and $D^{N}$, the Levi-Civita connection, the Clifford multiplication and the intrinsic Dirac operator on $N$. We can also define an extrinsic spinor bundle on $N$ by ${\bf{S}}:=\Sigma M_{|N}$. Then we recall the identification between these two spinor bundles ({\it cf} \cite{HMo}, \cite{Mo3} or \cite{Ba2} for instance):
\beqt\label{identifications} 
{\bf{S}}\equiv \left\lbrace 
\begin{array}{ll}
 \Sigma N&\text{if $n$ is even}  \\ 
\Sigma N\oplus\Sigma N & \text{if $n$ is odd}.
\end{array} \right. 
\eeqt
\noindent The interest of this identification is that we can use restrictions of ambient spinors to study the intrinsic Dirac operator of $N$. Indeed, we can define an  extrinsic connection $\nabla^{{\bf{S}}}$ and a Clifford multiplication $\gamma^{{\bf{S}}}$ on ${\bf{S}}$ by
\beqt\label{gaussspin} 
\nabla^{{\bf{S}}}=\nabla+\frac{1}{2}\gamma(\nu)\gamma(A),
\eeqt
\beqt
\gamma^{{\bf{S}}}=\gamma(\nu)\gamma,
\eeqt
\noindent where $\nu$ is the exterior normal unit vector field and $A$ the associated Weingarten operator. By the previous identification given by (\ref{identifications}), we can also identify connections and Clifford multiplications.
\beqt
\nabla^{{\bf{S}}}\equiv \left\lbrace \begin{array}{ll}
 \nabla^{N}& \text{if $n$ is even,} \\ 
\nabla^{N}\oplus \nabla^{N}& \text{if $n$ is odd,} 
\end{array} \right. 
\eeqt
\beqt
\gamma^{{\bf{S}}}\equiv \left\lbrace \begin{array}{ll}
 \gamma^{N}& \text{if $n$ is even,} \\ 
\gamma^{N}\oplus -\gamma^{N}& \text{if $n$ is odd.} 
\end{array} \right. 
\eeqt
\noindent Then, we can consider the following extrinsic Dirac operator ${\bf{D}}$ on $N$, acting on sections of ${{\bf{S}}}$ and given locally by
\beqt
{\bf{D}}=\sum_{i=1}^{n}\gamma^{{\bf{S}}}(e_i)\nabla^{{\bf{S}}}_{e_i},
\eeqt 
\noindent where $\base$ is an orthonormal local frame of
$TN$. Then, by (\ref{gaussspin}), we have
\beqt
{\bf{D}}=\frac{n}{2}{H}-\gamma(\nu)\sum_{i=1}^{n}\gamma(e_i)\nabla_{e_i},
\eeqt 
\noindent that is, for any $\psi\in\Gamma({\bf{S}})$
\beqt\label{Diracbord} 
{\bf{D}}\psi:=\frac{n}{2}{H}\psi-\gamma(\nu)D\psi-\nabla_{\nu}\psi.
\eeqt
\noindent Moreover, extrinsic and intrinsic Dirac operators are related by the following formula
\beqt
{\bf{D}}\equiv \left\lbrace \begin{array}{ll}
 D^{N}& \text{if $n$ is even,} \\ 
D^{N}\oplus -D^{N}& \text{if $n$ is odd.} 
\end{array} \right. 
\eeqt
\noindent It is easy to see from the definition that
\beqt\label{commute} 
\nabla^{{\bf{S}}}\gamma(\nu)=\gamma(\nu)\nabla^{{\bf{S}}}\quad\text{and}\quad{\bf{D}}\gamma(\nu)=-\gamma(\nu){\bf{D}}.
\eeqt
\noindent This implies that if $N$ is compact, then the spectrum of $ {\bf{D}}$ is symmetric and
\beqt
\sspec\,{{\bf{D}}} \equiv \left\lbrace \begin{array}{ll}
\sspec\left(D^{N}\right) & \text{if $n$ is even,}\\\\
\sspec\left(D^{N}\right)\cup-\sspec\left(D^{N}\right) & \text{if $n$ is odd.}
\end{array} \right. 
\eeqt
\subsection{3-dimensional homogeneous manifolds with\\
4-dimensional isometry group}
\label{sect21}
In this section, we give a description of the 3-dimensional homogeneous manifolds with 4-dimensional isometry group. Such a manifold is a Riemannian fibration over a simply connected 2-dimensional manifold with constant curvature $\kappa$ and such that the fibers are geodesic. We denote by $\tau$ the bundle curvature, which measures the defect of the fibration to be a Riemannian product. When $\tau$ vanishes, we get a product manifold $\M^2(\kappa)\times\R$. Here, we describe 3-homogeneous manifolds with 4-dimensional isometry group and  $\tau\neq0$. These manifolds are of three types: they have the isometry group of the Berger spheres if $\kappa>0$, of the Heisenberg group $Nil_3$ if $\kappa=0$ or of $\widetilde{PSL_2(\R)}$ if $\kappa<0$ (see Table 1).\\\\
\begin{center}
\begin{tabular}{|c|c|c|c|}
\hline
&$\kappa>0$&$\kappa=0$&$\kappa<0$\\
\hline
$\tau=0$&$\Ss^2(\kappa)\times\R$&&$\HH^2\times\R$\\
\hline
$\tau\neq0$&$(\Ss^3,g_{Berger})$&$Nil_3$&$\widetilde{PSL_2(\R)}$\\
\hline
\end{tabular}\\
\vspace{0.2cm}
Table 1: \textit{Classification of }$\Ekt$
\end{center}
\vspace{0.3cm}
In the sequel, we denote these homoegenous manifolds by $\Ekt$. We begin by giving a precise description of $\Ekt$. For further details, one can refer to \cite{Dan2} or \cite{Sco}.
\subsubsection{Canonical frame}
\label{sect211}
Let $\Ekt$ be a 3-dimensional homogeneous manifold with 4-dimensional isometry group. Assume that $\tau\neq0$, {\it i.e.}, $\Ekt$ is not a  product manifold $\mathbb{M}^2(\kappa)\times\R$. As we said, $\Ekt$ is a Riemannian fibration over a simply connected 2-dimensional manifold with constant curvature $\kappa$ and such that the fibers are geodesic. Now, let $\xi$ be a unitary vector field tangent to the fibers. We call it the vertical vector field. This vector field is a Killing vector field (corresponding to translations along the fibers).\\
\indent
We denote respectively by $\nablab$ and $\Rbar$ the Riemannian connection and the curvature tensor of $\Ekt$. The manifold $\Ekt$ admits a local direct orthonormal frame $\{e_1,e_2,e_3\}$ with 
$$e_3=\xi$$
and such that the Christoffel symbols $\overline{\Gamma}_{ij}^k=\left\langle \nablab_{e_i}e_j,e_k\right\rangle$ are
\beqt\label{christoffel}
\left\lbrace  
\begin{array}{l}
\overline{\Gamma}_{12}^3=\overline{\Gamma}_{23}^1=-\overline{\Gamma}_{21}^3=-\overline{\Gamma}_{13}^2=\tau,\\ \\
\overline{\Gamma}_{32}^1=-\overline{\Gamma}_{31}^2=\tau-\sigma, \\ \\
\overline{\Gamma}_{ii}^i=\overline{\Gamma}_{ij}^i=\overline{\Gamma}_{ji}^i=\overline{\Gamma}_{ii}^j=0,\quad\forall\,i,j\in\{1,2,3\},
\end{array}
\right. 
\eeqt
where $\sigma=\dfrac{\kappa}{2\tau}$. Then we have
$$[e_1,e_2]=2\tau e_3,\quad [e_2,e_3]=\sigma e_1,\quad [e_3,e_1]=\sigma e_2.$$
We will call $\{e_1,e_2,e_3\}$ the canonical frame of $\Ekt$. From (\ref{christoffel}), we see easily that for any vector field $X$,
\beqt
\nablab_Xe_3=\tau X\wedge e_3,
\eeqt
where $\wedge$ is the vector product in $\Ekt$, that is, for any $X,Y,Z\in\Gamma(T\,M)$,
$$\left\langle X\wedge Y,Z\right\rangle =\ddet_{\{e_1,e_2,e_3\}}(X,Y,Z).$$

\subsubsection{Manifolds with the isometry group of the Berger spheres}
\label{sect212}
Here, we are in the case $\kappa>0$ and $\tau\neq0$. These manifolds are fibrations over round 2-spheres. They are obtained by deforming the metric of a round 3-spheres in a way preserving the Hopf fibration and modifying the length of fibers.\\
\indent
We can see $\R^3$ endowed with the following metric 
$$\diff s^2=\lambda^2(\diff x^2+\diff y^2)+\big(\tau(y\diff x-x\diff y)+\diff z\big)^2$$
with
$$\lambda(x,y,z)=\frac{1}{1+\frac{\kappa}{4}(x^2+y^2)}$$
as the universal cover of such a homogeneous manifold ($\kappa>0$ et $\tau\neq0$) minus the fiber corresponding to the point $\infty\in\Ss^2$. In these coordinates, the fibers are $\left\lbrace x=x_0\ ,\ y=y_0\right\rbrace $. Then, the canonical frame is $\left\lbrace e_1,e_2,e_3\right\rbrace$ defined as follows.
$$\left\lbrace \begin{array}{l}
 e_1=\lambda^{-1}\big(\cos(\sigma z)\partial_x+\sin(\sigma z)\partial_y\big)+\tau\big(x\sin(\sigma z)- y\cos(\sigma z)\big)\partial_z,\\
 e_2=\lambda^{-1}\big(-\sin(\sigma z)\partial_x+\cos(\sigma z)\partial_y\big)+\tau\big(x\cos(\sigma z)+ y\sin(\sigma z)\big)\partial_z,\\ 
e_3=\xi,
\end{array}
\right. $$
with
$$\sigma=\frac{\kappa}{2\tau}.$$
Thus, this frame satisfies
$$[e_1,e_2]=2\tau e_3,\qquad [e_2,e_3]=\frac{\kappa}{2\tau} e_1,\qquad [e_3,e_1]=\frac{\kappa}{2\tau} e_2,$$
and is defined on the open set $\Ekt'$ which is $\Ekt$ minus the fiber corresponding to the point at the infinity of $\Ss^2(\kappa)$.
\begin{rem}
In the strict sense, the Berger spheres are $\Ekt$ with $\kappa=4$.
\end{rem}

\subsubsection{Manifolds with the isometry group of the Heisenberg space $Nil_3$}
\label{sect213}
They occur when $\kappa=0$ and $\tau\neq0$. They are fibrations over the Euclidean space $\R^2$ .\\
Recall that the Heisenberg group is the Lie group $Nil_3$ defined by
$$Nil_3=\left\lbrace \left( 
\begin{array}{ccc}
1 & a & b \\ 
0 & 1 &  c\\ 
0 & 0 & 1
\end{array}
\right) \ ;\ (a,b,c)\in\R^3\right\rbrace $$
endowed with a left invariant metric. Using exponential coordinates allows us to see the Heisenberg group as $\R^3$ with the following metric.
$$\diff s^2=\diff x^2+\diff y^2+\big(\tau(y\diff x-x\diff y)+\diff z\big)^2.$$
In these coordinates, the fibers are given by $\left\lbrace x=x_0\ ,\ y=y_0\right\rbrace $ and the canonical frame $\left\lbrace e_1,e_2,e_3\right\rbrace$ is defined by
$$\left\lbrace \begin{array}{l}
e_1=\partial_x-\tau y\partial_z, \\ 
e_2=\partial_y+\tau x\partial_z,\\ 
e_3=\partial_z.
\end{array}\right. $$
Moreover this frame satisfies
$$[e_1,e_2]=2\tau e_3,\qquad[e_2,e_3]=0,\qquad[e_3,e_1]=0.$$

\subsubsection{Manifolds with the isometry group of $\widetilde{PSL_2(\R)}$}
\label{sect214}
Such manifolds are fibrations over hyperbolic plans. They occur when $\kappa<0$ and $\tau\neq0$. We can take the Poincar\'e disk model for the hyperbolic plane of constant curvature $\kappa<0$. The manifold $\mathbb{D}^2\left( \frac{2}{\sqrt{-\kappa}}\right)\times\R$, endowed with the metric
$$\diff s^2=\lambda^2(\diff x^2+\diff y^2)+\big(\tau(y\diff x-x\diff y)+\diff z\big)^2,$$
with
$$\lambda(x,y,z)=\frac{1}{1+\frac{\kappa}{4}(x^2+y^2)},$$
and where $\mathbb{D}^2(\rho)$ is the open disk of radius $\rho$ in $\R^2$, is a homogeneous manifold of bundle curvature $\tau$ and base curvature $\kappa$. The fibers are $\left\lbrace x=x_0\ ,\ y=y_0\right\rbrace $. Then, the canonical frame is $\left\lbrace e_1,e_2,e_3\right\rbrace$ defined as follows
$$\left\lbrace \begin{array}{l}
 e_1=\lambda^{-1}\big(\cos(\sigma z)\partial_x+\sin(\sigma z)\partial_y\big)+\tau\big(x\sin(\sigma z)- y\cos(\sigma z)\big)\partial_z,\\
 e_2=\lambda^{-1}\big(-\sin(\sigma z)\partial_x+\cos(\sigma z)\partial_y\big)+\tau\big(x\cos(\sigma z)+ y\sin(\sigma z)\big)\partial_z,\\ 
e_3=\xi,
\end{array}
\right. $$
where
$$\sigma=\frac{\kappa}{2\tau},$$
satisfies
$$[e_1,e_2]=2\tau e_3,\qquad [e_2,e_3]=\frac{\kappa}{2\tau} e_1,\qquad [e_3,e_1]=\frac{\kappa}{2\tau} e_2.$$

\subsubsection{The curvature tensor of $\Ekt$}
\label{sect215}
From the description of the homogeneous spaces $\Ekt$ that we gave, we can explictly compute the curvature tensor. In the three cases, there exists a canonical frame $\{e_1,e_2,e_3\}$ such that the Christoffel symbols are given by (\ref{christoffel}). From these relations, we can see that the curvature operator $\Rbar$, acting on 2-forms, defined by
 $$\left\langle \overline{R}(X\wedge Y),Z\wedge W\right\rangle :=\left\langle \overline{R}(X,Y)Z,W\right\rangle $$
 is diagonal in the frame $\{e_2\wedge e_3,e_3\wedge e_1,e_1\wedge e_2\}$ and the corresponding matrix is
 \beqt\label{diagonal}
 \overline{R}={\rm diag}(\tau^2,\tau^2,\kappa-3\tau^2).
 \eeqt
Because of the symmetry of the curvature tensor, it is completely determined by (\ref{diagonal}). Namely, we have
\begin{prop}\label{courbure}
For $X,Y,Z,W\in\Chi(\Ekt)$, we have
$$\left\langle \overline{R}(X,Y)Z,W\right\rangle =(\kappa-3\tau^2)\left\langle R_0(X,Y)Z,W\right\rangle+(\kappa-4\tau^2)\left\langle  R_1(\xi;X,Y)Z,W\right\rangle ,$$
where $R_0$ and $R_1$ are defined by
$$R_0(X,Y)Z=\left\langle X,Z\right\rangle Y-\left\langle Y,Z\right\rangle X,$$
and
\beQ
R_1(V;X,Y)Z&=&\left\langle Y,V\right\rangle \left\langle Z,V\right\rangle X+\left\langle Y,Z\right\rangle \left\langle X,V\right\rangle V\\
&&-\left\langle X,Z\right\rangle \left\langle Y,V\right\rangle V-\left\langle X,V\right\rangle \left\langle Z,V\right\rangle Y
\eeQ
\end{prop}
\noindent
{\it Proof:} Any vector field $X$ can be decomposed in vertical and horizontal parts, $X=\widetilde{X}+x\xi$, where $x=\left\langle X,\xi\right\rangle $. By the same way, we decompose $Y$, $Z$ and $W$. Then, we get a sum of 16 terms for $\left\langle \overline{R}(X,Y)Z,W\right\rangle$. Using the symmerty of the curvature tensor, many terms vanish. Namely, when the vertical vector $\xi$ appears three times,  four times, or two times, but in position 1 and 2, or 3 and 4. Finally, by (\ref{diagonal}), when $\xi$ appears one time, the term also vanishes. Then, we get
\beQ
\left\langle \overline{R}(X,Y)Z,W\right\rangle&=&\left\langle \overline{R}(\widetilde{X},\widetilde{Y})\widetilde{Z},\widetilde{W}\right\rangle\\ 
&&+yw\left\langle \overline{R}(\widetilde{X},\xi)\widetilde{Z},\xi\right\rangle+yz\left\langle \overline{R}(\widetilde{X},\xi)\xi,\widetilde{W}\right\rangle\\
&&+xw\left\langle \overline{R}(\xi,\widetilde{Y})\widetilde{Z},\xi\right\rangle+xz\left\langle \overline{R}(\xi,\widetilde{Y})\xi,\widetilde{W}\right\rangle\\
&=&(\kappa-3\tau^2)\bigg(\left\langle \widetilde{X},\widetilde{Z}\right\rangle\left\langle \widetilde{Y},\widetilde{W} \right\rangle -\left\langle \widetilde{X},\widetilde{W}\right\rangle\left\langle \widetilde{Y},\widetilde{Z} \right\rangle\bigg)\\
&&+\tau^2\bigg(yw\left\langle \widetilde{X},\widetilde{Z}\right\rangle-yz\left\langle \widetilde{X},\widetilde{W}\right\rangle-xw\left\langle \widetilde{Y},\widetilde{Z}\right\rangle+xz\left\langle \widetilde{Y},\widetilde{W}\right\rangle\bigg)\\
&=&(\kappa-3\tau^2)\big(\left\langle X,Z\right\rangle \left\langle Y,W\right\rangle- \left\langle X,W\right\rangle \left\langle Y,Z\right\rangle \big)\\
&&-(\kappa-4\tau^2)\big(\left\langle X,Z\right\rangle \left\langle Y,\xi\right\rangle \left\langle W,\xi\right\rangle +\left\langle Y,W\right\rangle \left\langle X,\xi\right\rangle \left\langle Z,\xi\right\rangle \\
&&-\left\langle X,W\right\rangle \left\langle Y,\xi\right\rangle \left\langle Z,\xi\right\rangle -\left\langle Y,Z\right\rangle \left\langle X,\xi\right\rangle \left\langle W,\xi\right\rangle \big),
\eeQ
which achieves the proof.
$\finpreuve$
\section{Special spinor fields on the surfaces of $\Ekt$}
\subsection{Generalized Killing spinor for $\M^2(\kappa)\times\R$}
\subsubsection{Spinor bundle of $\MR$}
\label{sect23}
First, in every dimension, since $\M^n$ is spin, then
$\M^n\times\R$ is also spin and the spin structures
of $\M^n$ and $\M^n\times\R$ are in a one-to-one correspondence. Indeed,
If we consider a spin structure on $\M^n\times\R$, it is well-known that it induces a spin structure on $\M^n$. Conversely, if we consider a spin structure on $\M^n$, then we can lift it to a spin structure on $\M^n\times\R$, and we obtain a
$\widetilde{Gl}^+(n,\R)$-principal bundle over $\M^n\times\R$. Then we extend the structural groupe via the following embedding
$$\widetilde{Gl}^+(n,\R)\hookrightarrow\widetilde{Gl}^+(n+1,\R),$$
obtained by this standard embedding
$$Gl^+(n,\R)  \hookrightarrow  Gl^+(n+1,\R)$$
given by
$$a\longmapsto\left(\begin{array}{cc}
1 & 0\\
0 & a
\end{array}\right)$$
Finally, by restriction of this spin structure on $\M^n$, we get back the initial spin structure ({\it cf} \cite{BGM} for more details).\\
\indent
From now, we will focus on the case $n=2$. We will explain explictely the spinor bundle of $\M^2\times\R$. We have seen in Section \ref{sect22} that
$$\Sigma(\M^2\times\R)_{|\M^2}\cong\Sigma\M^2.$$
So, if $\varphi\in\Gamma( \Sigma(\M^2\times\R))$, then for any
$t\in\R$:
$$\varphi(.,t)\in\Gamma(\Sigma\M^2).$$
So we have the following vector-bundle isomorphism
\function{\beta}{\Sigma(\M^2\times\R)}{\Sigma(\M^2)\oplus\R}{\varphi}{(\widetilde{\varphi,0),}}
where $\oplus$ is a direct sum of vector-bundle and $\R$ the vector bundle
on $\R$ with fibre $\{0\}$. It means that the fibre $\Sigma_{(x,t)}(\M^2\times\R)$ over every point $(x,t)\in\M^2\times\R$ is the fibre $\Sigma_x\M^2$ of the spinor bundle $\M^2$.\\
Then a section of $\Sigma(\M^2\times\R)$ is a $\Cinf$-map
\function{\varphi}{\R}{\Gamma(\Sigma(\M^2))}{t}{\varphi_t.}
\subsubsection{Restriction to a surface}
\label{sec24}
In the sequel, we will consider some particular and interesting sections, precisely the sections which do not depend on $t$.\\
\indent
We have the spinorial Gauss formula,  
$$\nablab_X\varphi=\nabla_X\varphi+\frac{1}{2}\gamma(AX)\gamma(\nu)\varphi,$$
where $\nablab$ is the spinorial connection on $\M^2\times\R$, the spinorial connection on  $\M^2$ is $\nabla$, the Clifford multiplication on $\M^2\times\R$ is $\gamma$ and $A$ the Weingarten operator of the immersion of $\M^2$ into $\M^2\times\R$. Since $\M^2$ is totally geodesic in the product $\M^2\times\R$, we get by taking $\varphi_t=\varphi_0$ a Killing spinor on $\M^2$, {\it i.e.} $\nabla_X\varphi_0=\eta\gamma^{\M^2}(X)\varphi_0$:
\beQ
\nablab_X\varphi&=&\nabla_X\varphi_0\\
&=&\eta\gamma^{\M^2}(X)\varphi_0\\
&=&\eta\gamma(X)\gamma\left( \dt\right) \varphi_0\\
&=&\eta\gamma(X)\gamma(\dt)\varphi
\eeQ
On the other hand, the complexe volume form $\omega_{\C}=-e_1\cdot e_2\cdot\dt$ acts as identity, so we have
\beQ
\gamma(e_1)\gamma\left( \dt\right) \varphi&=&\gamma(e_2)\gamma(e_1)\gamma(e_2)\gamma\left( \dt\right) \varphi\\
&=&-\gamma(e_2)\varphi.
\eeQ
Similarly,
\beQ
\gamma(e_2)\gamma\left( \dt\right) \varphi&=&-\gamma(e_1)\gamma(e_1)\gamma(e_2)\gamma\left( \dt\right) \varphi\\
&=&\gamma(e_1)\varphi.
\eeQ
So we deduce that
\beqt\label{spineurM2R}
\left\lbrace \begin{array}{l}
\nablab_{e_1}\varphi=-\eta\gamma(e_2)\varphi, \\ \\
\nablab_{e_2}\varphi=\eta\gamma(e_1)\varphi,\\ \\
\nablab_{\dt}\varphi=0.
\end{array}\right. 
\eeqt
These particular spinor fields are the analogue of Killing spinor field for $\Ss^3$ or $\HH^3$. They are the most special spinors we can find on the product $\MR$. \\
\indent
Now, let $(N,\left\langle \cdot,\cdot\right\rangle )$ be a surface of $\M^2\times\R$, oriented by $\nu$. Since $(N,\left\langle \cdot,\cdot\right\rangle )$ is oriented, it could be equiped with a spin structure induce from the spin structure of $\M^2\times\R$. Moreover, as we saw, we have the following identification between the spinor bundles
$$\Sigma\,(\M^2\times\R)_{|N}\cong\Sigma\,N,$$
and the spinorial Gauss forumla (\ref{gaussspin}) gives the relation between the spinorial connctions of $N$ and $\M^2\times\R$. For any $X\in\Chi(N)$ and any $\psi\in\Gamma(\Sigma\,(\M^2\times\R))$, we have
\beQ
\big(\overline{\nabla}_X\psi\big)_{\big|N}&=&\nabla_X\big(\psi_{\big|N}\big)+\frac{1}{2}\gamma^N(AX)\psi_{\big|N}\\
&=&\nabla_X\big(\psi_{\big|N}\big)+\frac{1}{2}\gamma(AX)\gamma(\nu)\psi_{\big|N},
\eeQ
where $\gamma$ is the Clifford multiplication in $\M^2\times\R$, $\gamma^N$ the Clifford multiplication in $N$ and $A$ the Weingarten operator of the immersion of $N$ into $\M^2\times\R$. If we use this forumla for the special spinor field on $ \M^2\times\R$ given by (\ref{spineurM2R}), we get
\beQ
\nabla_X\varphi&=&\nablab_X\varphi-\frac{1}{2}\gamma^N(AX)\varphi\\
&=&\eta\gamma(X_t)\gamma\left( \dt\right) \varphi-\frac{1}{2}\gamma^N(AX)\varphi,
\eeQ
where $X_t$ is the part of $X$ tangent to $\M^2$, that is,
\beQ
X_t&=&X-\left\langle X,\dt\right\rangle\dt\\
&=&X-\left\langle  X,T\right\rangle T-f\left\langle X,T\right\rangle\nu. 
\eeQ 
So, we deduce that
\beQ
\nabla_X\varphi&=&\eta\gamma(X)\gamma\left( \dt\right)\varphi-\eta\left\langle X,\dt\right\rangle \gamma\left( \dt\right)\gamma\left( \dt\right)\varphi-\frac{1}{2}\gamma^N(AX)\varphi\\
&=&\eta\gamma(X)\gamma\left( \dt\right)\varphi+\eta\left\langle X,\dt\right\rangle \varphi-\frac{1}{2}\gamma^N(AX)\varphi\\
&=&\eta\gamma(X)\gamma(T)\varphi+\eta f\gamma(X)\gamma(\nu)\varphi+\eta\left\langle X,T\right\rangle \varphi-\frac{1}{2}\gamma^N(AX)\varphi.
\eeQ
On the other hand, if we denote by $\omega=e_1\underset{N}{\cdot}e_2$ the real volume element on $N$, we have the following relations
$$\left\lbrace \begin{array}{l}
 \gamma(X)=-\gamma^N(X)\gamma^N(\omega),\\ 
\gamma(\nu)=\gamma^N(\omega).
\end{array}\right. $$
By using these two identities, the fact that $\omega^2=-1$ and that $\omega$ anti-commuts with vector fields tangent to $N$, we get
\beQ
\nabla_X\varphi&=&\eta\gamma^N(X)\gamma^N(\omega)\gamma^N(T)\gamma^N(\omega)\varphi+\eta\left\langle X,T\right\rangle \varphi-\eta f\gamma^N(X)\gamma^N(\omega)\gamma^N(\omega)\varphi\\
&&-\frac{1}{2}\gamma^N(AX)\varphi
\eeQ
Now, we can rewrite this equation as follows
\beq\label{killinggeneral}
\nabla_X\varphi&=&
\eta X\cdot T\cdot\varphi+\eta fX\cdot\varphi+\eta\left\langle X,T\right\rangle \varphi-\frac{1}{2}AX\cdot\varphi.
\eeq
where ``$\cdot$'' stands for the Clifford multiplication on $N$.
\begin{defn}
A spinor field which satifies the equation (\ref{killinggeneral}) is called a generalized Killing spinor.
\end{defn}
These spinor fields satisfy the following property
\begin{prop}\label{normeconstante1}
\be[i)]
\item If $\eta=\frac{1}{2}$, then the norm of a generalized Killing spinor is constant.
\item If $\eta=\frac{i}{2}$, then the norm of a generalized Killing spinor satisfies for any $X\in\Chi(N)$:
$$X|\varphi|^2= \pre\left\langle  iX\cdot T\cdot\varphi+ifX\cdot\varphi,\varphi\right\rangle. $$
\ee
\end{prop}
\noindent
{\it Proof:} We need to calculate $X|\varphi|^2$ for $X\in\Chi(N)$. We have,
$$X|\varphi|^2=2\pre\left\langle \nabla_X\varphi,\varphi\right\rangle .$$
We replace $\nabla_X\varphi$ by the expression given by (\ref{killinggeneral}), and we use the following lemma
 \begin{lem}\label{pdtscalaireformesreelles}
 Let $\psi$ be a spinor field and $\beta$ a real 1-form or 2-form. Then
 $$\pre\left\langle \beta\cdot\psi,\psi\right\rangle =0.$$
 \end{lem}
\noindent
 We deduce easily form this lemma that
 $$\pre\left\langle A(X)\cdot\varphi,\varphi\right\rangle =0,$$
 and
 $$\pre\left\langle fX\cdot\varphi,\varphi\right\rangle =0.$$
 By this lemma again, we see that
 $$\pre\left\langle X\cdot T\cdot\varphi,\varphi\right\rangle +\pre\left\langle \left\langle X,T\right\rangle\varphi,\varphi \right\rangle=0.$$ 
If $\eta=\frac{i}{2}$, the first three terms are non zero but an simple calculation yields to the result. $\finpreuve$
\begin{rem}
In the case $\eta=\frac{i}{2}$, the norm of $\varphi$ is not constant. Nevertheless, we can show that $\varphi$ never vanishes.
\end{rem}
\subsection{Generalized Killing spinors on $\Ekt$, $\tau\neq0$}
\subsubsection{Special spinor fields on $\Ekt$, $\tau\neq0$}
\label{sect23}
We will give here the expression of special spinor fields on $\Ekt$ which will play the role of parallel or Killing spinors that we have in space forms. For this, let us endow $\Ekt$ with its trivial spinorial structure. We can consider a constant section $\varphi$ of the spinor bundle. From the local expression of the spinorial connection and of the Christoffel symbols $\Gamma_{ij}^k$ given in Section \ref{sect211}, we can compute $\overline{\nabla}_X\varphi$. We recall that the complex volume element $\omega_{3}^{\C}=-e_1\cdot e_2\cdot\xi$ acts as identity, which implies
\beQ
\gamma(e_1)\gamma(e_2)\varphi&=&-\gamma(e_1)\gamma(e_2)\gamma(e_1)\gamma(e_2)\gamma(\xi)\varphi\\
&=&\gamma(\xi)\varphi.
\eeQ
\noindent
Similarly, we have
$$\gamma(e_2)\gamma(\xi)\varphi=\gamma(e_1)\varphi \quad\text{and}\quad \gamma(\xi)\gamma(e_1)\varphi=\gamma(e_2)\varphi.$$
The expression of the spinorial connection is for any $X\in\Chi(\Ekt)$
\beQ
\overline{\nabla}_X\varphi=X(\varphi)+\frac{1}{4}\sum_{i,j=1}^3\left\langle \nabla_Xe_i,e_j\right\rangle \gamma(e_i)\gamma(e_j)\varphi, 
\eeQ
So, we deduce that
\beQ
\overline{\nabla}_{e_1}\varphi&=&\frac{1}{4}\left( \Gamma_{12}^3\gamma(e_2)\gamma(\xi)\varphi+\Gamma_{13}^2\gamma(\xi)\gamma(e_2)\varphi\right)\\
&=&\frac{1}{2}\tau \gamma(e_2)\gamma(\xi)\varphi\\
&=&\frac{1}{2}\tau \gamma(e_1)\varphi.
\eeQ
By the same way, we have
$$\overline{\nabla}_{e_2}\varphi=\frac{1}{2}\tau \gamma(e_2)\varphi.$$
Then, 
\beQ
\overline{\nabla}_{\xi}\varphi&=&\frac{1}{4}\left( \Gamma_{31}^2\gamma(e_1)\gamma(e_2)\varphi+\Gamma_{32}^1\gamma(e_2)\gamma(e_1)\varphi\right)\\
&=&\frac{1}{2}\left(\frac{\kappa}{2\tau}-\tau\right)  \gamma(e_1)\gamma(e_2)\varphi\\
&=&\frac{1}{2}\left(\frac{\kappa}{2\tau}-\tau\right)\gamma(\xi)\varphi.
\eeQ
Finally, we have a special spinor field $\varphi$ satisfying
\beqt\label{presquekilling}
\left\lbrace \begin{array}{l}
 \overline{\nabla}_{e_1}\varphi=\dfrac{1}{2}\tau \gamma(e_1)\varphi,\\ \\
 \overline{\nabla}_{e_2}\varphi=\dfrac{1}{2}\tau \gamma(e_2)\varphi,\\ \\
\overline{\nabla}_{\xi}\varphi=\dfrac{1}{2}\left(\dfrac{\kappa}{2\tau}-\tau\right)\gamma(\xi)\varphi.
\end{array}\right.
\eeqt
\subsubsection{Restriction to a surface}
\label{sec24}
In this section, we give the restriction of these spinors to a surface of $\Ekt$. Let $(N,\left\langle\cdot,\cdot\right\rangle)$ be an orientable surface into $\Ekt$, oriented by $\nu$. Since $(N,\left\langle\cdot,\cdot\right\rangle)$ is oriented, it is equiped with a spinorial structure induced from the spinorial structure of $\Ekt$. Moreover, we recall the identification between spinor bundles
$$\Sigma\,\Ekt_{\big|N}\cong\Sigma\,N.$$
We also recall  the spinorial Gauss formula:
\beQ
\big(\overline{\nabla}_X\psi\big)_{\big|N}&=&\nabla_X\big(\psi_{\big|N}\big)+\frac{1}{2}\gamma^N(AX)\psi\\
&=&\nabla_X\big(\psi_{\big|N}\big)+\frac{1}{2}\gamma(AX)\gamma(\nu)\psi.
\eeQ
For the special spinor field satisying (\ref{presquekilling}), the Gauss formula yields to
$$\nabla_X\big(\varphi_{\big|N}\big)=\frac{\tau}{2}\gamma(\widetilde{X})\varphi_{\big|N}+\frac{1}{2}\left( \frac{\kappa}{2\tau}-\tau\right)x\gamma(\xi)\varphi_{\big|N} -\frac{1}{2}\gamma^N(AX)\varphi_{\big|N},$$
where $\widetilde{X}$ et $x$ are defined as in Section \ref{sect215} by
$$X=\widetilde{X}+x\xi,$$
that is, $x=\left\langle X,T\right\rangle$ and
\beQ
\widetilde{X}&=&X-\left\langle X,T\right\rangle\xi\\
&=&X-\left\langle X,T\right\rangle T -f \left\langle X,T\right\rangle\nu.
\eeQ
Finally, we get
\beQ
\nabla_X\varphi&=&\frac{\tau}{2}\gamma(X)\varphi-\frac{\tau}{2}\left\langle X,T\right\rangle\gamma(T)\varphi-\frac{\tau}{2}f\gamma(\nu)\varphi\\
&&+\frac{1}{2}\left( \frac{\kappa}{2\tau}-\tau\right)\left\langle X,T\right\rangle\gamma(T)\varphi+\frac{1}{2}\left( \frac{\kappa}{2\tau}-\tau\right)f\left\langle X,T\right\rangle\gamma(\nu)\varphi-\frac{1}{2}\gamma^N(AX)\varphi\\\\
&=&\frac{\tau}{2}\gamma(X)\varphi-\frac{\alpha}{2}\left\langle X,T\right\rangle\gamma(T)\varphi-\frac{\alpha}{2}f\left\langle X,T\right\rangle\gamma(\nu)\varphi-\frac{1}{2}\gamma(AX)\gamma(\nu)\varphi,
\eeQ
where we set $\alpha=2\tau-\frac{\kappa}{2\tau}$.
On the other hand, if we denote by $\omega=e_1\underset{N}{\cdot}e_2$ the real volume element on $N$, we have the following identities
$$\left\lbrace \begin{array}{l}
 \gamma(X)=-\gamma^N(X)\gamma^N(\omega),\\ 
\gamma(\nu)=\gamma^N(\omega).
\end{array}\right. $$
Then, using the fact that $\omega^2=-1$, we obtain
\beqt\label{killinggeneralise}
\nabla_X\varphi=-\frac{\tau}{2}X\cdot\omega\cdot\varphi+\frac{\alpha}{2}\left\langle X,T\right\rangle T\cdot\omega\cdot\varphi-\frac{\alpha}{2}f \left\langle X,T\right\rangle\omega\cdot\varphi-\frac{1}{2}AX\cdot\varphi,
\eeqt
where ``$\cdot$'' stands for the Clifford multiplication on $N$. \\
\noindent
We recall that the spinor bundle $\Sigma\,N$ decomposes in
$$\Sigma\,N=\Sigma^+\,N\oplus\Sigma^-\,N,$$
where $\Sigma^{\pm}\,N$ is the eigenspace for the eigenvalue $\pm1$ under the action of the complex volume element $\omega_2=i\omega$. In this decomposition, $\varphi$ can be written $\varphi=\varphi^++\varphi^-$, and we set $$\overline{\varphi}=\omega_2\cdot\varphi=\varphi^+-\varphi^-.$$
Then, the equation (\ref{killinggeneralise}) becomes
\beqt\label{killinggeneralise2}
\nabla_X\varphi=i\frac{\tau}{2}X\cdot\overline{\varphi}-i\frac{\alpha}{2}\left\langle X,T\right\rangle T\cdot\overline{\varphi}+i\frac{\alpha}{2}f\left\langle X,T\right\rangle\overline{\varphi}-\frac{1}{2}AX\cdot\varphi.
\eeqt
\begin{defn}
A spinor field satisfying (\ref{killinggeneralise2}) is called a {\it generalized Killing spinor}.
\end{defn}
\begin{prop}\label{normeconstante}
The norm of a generalized Killing spinor is constant.
\end{prop}
\noindent
{\it Proof:} It is sufficient to compute $X|\varphi|^2$ for any $X\in\Chi(N)$. Indeed,
$$X|\varphi|^2=2\pre\left\langle \nabla_X\varphi,\varphi\right\rangle .$$
We replace $\nabla_X\varphi$ by its expression given by (\ref{killinggeneralise}), and we use Lemma \ref{pdtscalaireformesreelles} to conclude.
 
$\finpreuve$
\section{Isometric immersions into $\M^2(\kappa)\times\R$}
\label{sect3}
\subsection{The results}
\label{sect34}
Here, we state the main results of the present paper which give a spinorial characterization of surfaces isometrically immersed into the product spaces $\Ss^2\times\R$ and $\HH^2\times\R$. Our first theorem concerns the surfaces into $\Ss^2\times\R$.
\begin{thm}\label{thm1} 
Let $(N,\langle.,.\rangle)$ be a connected, oriented and simply connectd Riemannian surface. Let $T$ be a vector field, $f$ and $H$ two real  functions on $N$ satisfying
$$\left\lbrace  \begin{array}{l}
f^2+||T||^2=1, \\ 
\langle\nabla_XT,Y\rangle=\langle\nabla_YT,X\rangle,\quad\forall X,Y\in\chi(N), \\ 
2Hf=\ddiv(T).
\end{array}\right. $$
Let $A\in{\rm{End}}(TN)$. The following three data are equivalent:
\begin{enumerate}[i)]
\item an isometric immersion $F$ from $N$ into $\Ss^2\times\R$ of mean curvature $H$ such that the Weingarten operator related to the normal $\nu$ is given by
$$dF\circ A\circ dF^{-1}$$ 
and such that
$$\dt=dF(T)+f\nu.$$
\item a spinor field $\varphi$ satisfying
$$\nabla_X\varphi=\frac{1}{2} X\cdot T\cdot\varphi+\frac{1}{2} fX\cdot\varphi+\frac{1}{2}\left\langle X,T\right\rangle \varphi-\frac{1}{2}AX\cdot\varphi,$$
where $A$ satisfies
$$\nabla_XT=fAX,$$
and
$$df(X)=-\langle AX,T\rangle.$$
\item a spinor field $\varphi$ satisfying
$$D\varphi=(H-f)\varphi-\frac{1}{2} T\cdot\varphi,$$ of constant norm and such that
$$df=-2Q_{\varphi}(T)+fT.$$
\end{enumerate}
\end{thm}
Here, $Q_{\varphi}$ is the energy-momentum tensor associated with the spinor field $\varphi$ and defined by
$$\qphi(X,Y):=\frac{1}{2}\pre\left\langle X\cdot\nabla_Y\varphi+Y\cdot\nabla_X\varphi,\varphi/|\varphi|^2\right\rangle .$$
We have the analogue for surfaces into $\HH^2\times\R$
\begin{thm}\label{thm2} 
Let $(N,\langle.,.\rangle)$ be a connected, oriented and simply connectd Riemannian surface. Let $T$ be a vector field, $f$ and $H$ two real  functions on $N$ satisfying
$$\left\lbrace  \begin{array}{l}
f^2+||T||^2=1, \\ 
\langle\nabla_XT,Y\rangle=\langle\nabla_YT,X\rangle,\quad\forall X,Y\in\chi(N), \\ 
2Hf=\ddiv(T).
\end{array}\right. $$
Let $A\in{\rm{End}}(TN)$. The following three data are equivalent:
\begin{enumerate}[i)]
\item an isometric immersion $F$ from $N$ into $\HH^2\times\R$ of mean curvature $H$ such that the Weingarten operator related to the normal $\nu$ is given by
$$dF\circ A\circ dF^{-1}$$ 
and such that
$$\dt=dF(T)+f\nu.$$
\item a spinor field $\varphi$ satisfying
$$\nabla_X\varphi=\frac{i}{2} X\cdot T\cdot\varphi+\frac{i}{2} fX\cdot\varphi+\frac{i}{2}\left\langle X,T\right\rangle \varphi-\frac{1}{2}AX\cdot\varphi,$$
where $A$ satisfies
$$\nabla_XT=fAX,$$
and
$$df(X)=-\langle AX,T\rangle.$$
\item a spinor field $\varphi$ which never vanishes and satisfying
 $$D\varphi=(H_if)\varphi-\frac{i}{2} T\cdot\varphi,$$ 
 
$$X|\varphi|^2=\pre\big( iX\cdot T\cdot\varphi+ifX\cdot\varphi\varphi\big),$$
and such that
$$df=-2Q_{\varphi}(T)+B(T),$$
where $B$ is defined as in Lemma \ref{lienAQphi}
\end{enumerate}
\end{thm}

\subsection{Compatibility equation in $\M^n\times\R$}
\label{sect21}
Before giving the proof of these two theorems, we need to recall some facts about how to obtain isometric immersions into product spaces.
\\
Let $N$ be an orientable hypersurface of $\MRn$ and $\nu$ its unit normal vector. Let $T$ be the projection of the vector $\dt$ on the tangent bundle $T\,N$. Moreover, we consider the function $f$ defined by:
$$f:=\Big<\nu,\Dt\Big>.$$
It is clear that
$$\Dt=T+f\nu.$$
Since $\dt$ is a unit vector field, we have:
$$||T||^2+f^2=1.$$
Let's compute the curvature tensor of $\MRn$ for tangent vectors to $N$.
\begin{prop}
For all $X,Y,Z,W\in\Gamma(TN)$, we have:
\begin{eqnarray*}
\big<\overline{R}(X,Y)Z,W\big>&=&\kappa\big(\langle X,Z\rangle \langle Y,W\rangle -\langle Y,Z\rangle \langle X,W\rangle \\
&&-\langle Y,T\rangle \langle W,T\rangle \langle X,Z\rangle -\langle X,T\rangle \langle Z,T\rangle \langle Y,W\rangle \\
&&+\langle X,T\rangle \langle W,T\rangle \langle Y,Z\rangle +\langle Y,T\rangle \langle Z,T\rangle \langle X,W\rangle \big),
\end{eqnarray*}
and
\begin{eqnarray*}
\big<\overline{R}(X,Y)\nu,Z\big>&=&\kappa f\big(\langle X,Z\rangle \langle Y,T\rangle -\langle Y,Z\rangle \langle X,T\rangle \big).
\end{eqnarray*}
\end{prop}
\pf Let $X$ be a vector field over $\MRn$. It can be written as $X(m,t)=(X_t(m),X_m(t))$ with for all $t\in\R$, $X_t$ is a vector field over $\M^n$ and for all $m\in\M^n$, $X_m$ is a vector field over $\R$. In other words, we have:
\begin{eqnarray*}
X(m,t)&=&X_t(m)+\Big<X,\Dt\Big>\Dt(t)\\
&=&X_t(m)+\Big<X,T\Big>\Dt(t).
\end{eqnarray*}
Since $\dt$ is parallel, we get, for all $X,Y,Z,W\in\Chi(\MRn)$:
\begin{eqnarray*}
\overline{R}(X,Y,Z,W)&=&\langle\overline{R}_{\mathbb{M}^n}(X_t,Y_t)Z_t,W_t\rangle\\
&=&\kappa\big(\langle X_t,Z_t\rangle \langle Y_t,W_t\rangle -\langle Y_t,Z_t\rangle \langle X_t,W_t\rangle \big).
\end{eqnarray*}
Using the fact that
$$X_t(m)=X-\Big<X,T\Big>\Dt(t),$$
we obtain the first part of the proposition.\\
The second identity is proved identically with
$$\nu_t(m)=\nu-f\Dt(t).$$
$\pfend$
\begin{rem}\label{dim2}
If the couple $(X,Y)$ is orthonormal, we get
$$ \overline{R}(X,Y,X,Y)\big>=\kappa(1-\langle Y,T\rangle ^2-\langle X,T\rangle ^2),$$
which gives for $n=2$
$$\overline{R}_{1212}=\kappa(1-||T||^2)=\kappa f^2.$$
\end{rem}
\noindent
The fact that $\Dt$ is parallel implies the two following identities
\begin{prop}
For $X\in\Chi(N)$, we have
$$
\nabla_XT=fAX,
$$
and
$$
df(X)=-\langle AX,T\rangle.
$$

\end{prop}
\noindent
{\it Proof:} We know that $\overline{\nabla}_X\dt=0$ and $\dt=T+f\nu$,
so
\begin{eqnarray*}
0&=&\overline{\nabla}_XT+df(X)\nu+f\overline{\nabla}_X\nu\\
&=&\nabla_XT+\langle AX,T\rangle \nu+df(X)\nu-fAX.
\end{eqnarray*}
Now, it is sufficient to consider the normal and tangential parts to obtain the above identities.
$\finpreuve$
\begin{defn}[Compatibility Equations]\label{comp}
We say that $(N,\langle.,.\rangle,A,T,f)$ satisfies the compatibility equations for $\M^n\times\R$ if and only if for any $X,Y,Z\in\Chi(N)$,

\begin{align}\label{gauss}
R(X,Y)Z=&\langle AX,Z\rangle AY-\langle AY,Z\rangle AX\\
&+\kappa\Big(\langle X,Z\rangle Y-\langle Y,Z\rangle X-\langle Y,T\rangle \langle X,Z\rangle T\nonumber\\
&-\langle X,T\rangle \langle Z,T\rangle Y+\langle X,T\rangle \langle Y,Z\rangle T
+\langle Y,T\rangle \langle Z,T\rangle X\Big),\nonumber
\end{align}

\beqt\label{codazzi}
\nabla_XAY-\nabla_YAX-A[X,Y]=\kappa f(\langle Y,T\rangle X-\langle X,T\rangle Y),
\eeqt

\beqt\label{cond1}
\nabla_XT=fAX,\ \text{and}
\eeqt
\beqt\label{cond2}
df(X)=-\langle AX,T\rangle.
\eeqt
\end{defn}
\begin{rem}
The relations (\ref{gauss}) and (\ref{codazzi}) are the Gauss and Codazzi equations for an isometric immersion into $\M^n\times\R$.
\end{rem}
\subsection{A necessary condition}
\label{sect31}
Here, we recall a result of B. Daniel (\cite{Dan}) which gives a necessary and sufficient condition for the existence of an isometric immersion of an oriented, simply connected surface $N$ into $\Ss^2\times\R$ or $\HH^2\times\R$.
\begin{thrm}[Daniel \cite{Dan}]\label{dan2}
Let $(N,\langle .,.\rangle )$ be an oriented, simply connected surface and $\nabla$ its Riemannian connection. Let $A$ be field of symmetric endomorphisms $A_y:T_yN\longrightarrow T_yN$,
$T$ a vector field on $N$ and $f$ a smooth function on $N$, such that $||T||^2+f^2=1$. If $(N,\langle .,.\rangle ,A,T,f)$ satisfies the compatibility equations for $\M^2\times\R$, then, there exists an isometric immersion
$$F:N\longrightarrow\M^2\times\R$$ so that the Weingarten operator of the immersion related to the normal $\nu$ is
$$dF\circ A\circ dF^{-1}$$
and such that
$$\dt=dF(T)+f\nu.$$
Moreover, this immersion is unique up to a global isometry of
$\M^2\times\R$ which preserves the orientation of $\R$.
\end{thrm}

\subsection{Generalized Killing spinors and compatibility equations}
\label{sect32}
First, we show that the existence of a generalized Killing spinor implies the Gauss and Codazzi equations. For this, let $(N,\left\langle\cdot,\cdot\right\rangle)$ be an oriented surface with a generalized Killing spinor, that is, a non-trivial spinor field solution of the equation (\ref{killinggeneral}). We will see that the integrability conditions for this equation are precisely the Gauss and Codazzi equations.
\begin{prop}
Let $(N,\left\langle \cdot,\cdot\right\rangle )$ be an oriented surface with a non-trivial solution of the equation (\ref{killinggeneral}) and such that the equations (\ref{cond1}) and (\ref{cond2}) are satisfied. Then, the Gauss and Codazzi equations for $\M^2\times\R$ are also satisfied.
\end{prop}
\noindent
{\it Proof:} The proof of this proposition is based on the computation of the spinorial curvature applied to the spinor field $\varphi$ solution of the equation (\ref{killinggeneral}), {\it i.e.}
$$\mathcal{R}(X,Y)\varphi=\nabla_X\nabla_Y\varphi-\nabla_Y\nabla_X\varphi-\nabla_{[X,Y]}\varphi,$$
for $X,Y\in\Chi(N)$. Using the expression given by (\ref{killinggeneral}), the equations (\ref{cond1}) and (\ref{cond2}), the fact that $\omega^2=-1$ and that $\omega$ anticommutes with vector fields tangent to $N$, we get
\beQ
\nabla_X\nabla_Y\varphi&=&\underbrace{\eta fY\cdot AX\cdot\varphi}_{\alpha_1(X,Y)}+\underbrace{\eta^2Y\cdot T\cdot X\cdot T\cdot\varphi}_{\alpha_2(X,Y)}+\underbrace{\eta^2fY\cdot T\cdot X\cdot\varphi}_{\alpha_3(X,Y)}\\
&&-\underbrace{\frac{\eta}{2}Y\cdot T\cdot AX\cdot \varphi}_{-\alpha_{4}(X,Y)}-\underbrace{\eta\left\langle AX,T\right\rangle Y\cdot\varphi}_{-\alpha_{5}(X,Y)}+\underbrace{\eta^2fY\cdot X\cdot T\cdot\varphi}_{\alpha_{6}(X,Y)}\\
&&+\underbrace{\eta^2\left\langle X,T\right\rangle Y\cdot T\cdot\varphi}_{\alpha_{7}(X,Y)}+\underbrace{\eta^2f^2Y\cdot X\cdot\varphi}_{\alpha_{8}(X,Y)}+\underbrace{\eta^2f\left\langle X,T\right\rangle Y\cdot\varphi}_{\alpha_{9}(X,Y)}\\
&&-\underbrace{\frac{\eta}{2}fY\cdot AX\cdot\varphi}_{-\alpha_{10}(X,Y)}+\underbrace{\eta f\left\langle Y,AX\right\rangle\varphi }_{\alpha_{11}(X,Y)}+\underbrace{\eta^2\left\langle Y,T\right\rangle X\cdot T\cdot\varphi}_{\alpha_{12}(X,Y)}\\
&&+\underbrace{\eta^2f\left\langle Y,T\right\rangle X\cdot\varphi}_{\alpha_{13}(X,Y)}+\underbrace{\eta^2\left\langle X,T\right\rangle\left\langle Y,T\right\rangle \varphi}_{\alpha_{14}(X,Y)}-\underbrace{\frac{\eta}{2}\left\langle Y,T\right\rangle AX\cdot\varphi }_{-\alpha_{15}(X,Y)}\\
&&-\underbrace{\frac{1}{2}\nabla_X(AY)\cdot\varphi}_{-\alpha_{16}(X,Y)}-\underbrace{\frac{\eta}{2}AY\cdot X\cdot T\cdot\varphi}_{-\alpha_{17}(X,Y)}-\underbrace{\frac{\eta}{2}fAY\cdot X\cdot\varphi}_{-\alpha_{18}(X,Y)}\\
&&-\underbrace{\frac{\eta}{2}\left\langle X,T\right\rangle AY\cdot\varphi}_{-\alpha_{19}(X,Y)}+\underbrace{\frac{1}{4}AY\cdot AX\cdot\varphi}_{\alpha_{20}(X,Y)}+\underbrace{\eta\nabla_XY\cdot T\cdot\varphi}_{\alpha_{21}(X,Y)}\\
&&+\underbrace{\eta f\nabla_XY\cdot\varphi}_{\alpha_{22}(X,Y)}+\underbrace{\eta\left\langle \nabla_XY,T\right\rangle\varphi }_{\alpha_{23}(X,Y)}
\eeQ
That is, 
$$\nabla_X\nabla_Y\varphi=\sum_{i=1}^{23}\alpha_i(X,Y).$$
Obviously, by symmetry, we have
$$\nabla_Y\nabla_X\varphi=\sum_{i=1}^{23}\alpha_i(Y,X).$$
On the other hand, we have
\beQ
\nabla_{[X,Y]}\varphi&=&\underbrace{\eta[X,Y]\cdot T\cdot\varphi}_{\beta_1([X,Y])}+\underbrace{\eta f[X,Y]\cdot\varphi}_{\beta_2([X,Y])}\\
&&+\underbrace{\eta\left\langle  [X,Y],T\right\rangle\varphi}_{\beta_3([X,Y])}-\underbrace{\frac{1}{2}A[X,Y]\cdot\varphi}_{-\beta_4([X,Y])}.
\eeQ
Since the connection $\nabla$ is torsion-free, we get
$$\nabla_XY-\nabla_YX-[X,Y]=0,$$
which implies that
\beQ
\alpha_{21}(X,Y)-\alpha_{21}(Y,X)-\beta_1([X,Y])=0,\\
\alpha_{22}(X,Y)-\alpha_{22}(Y,X)-\beta_2([X,Y])=0,\\
\alpha_{23}(X,Y)-\alpha_{23}(Y,X)-\beta_3([X,Y])=0.
\eeQ
Moreover, by symmetry, we have
\beQ
\alpha_{11}(X,Y)-\alpha_{11}(Y,X)=0
\eeQ
and
\beQ
\alpha_{14}(X,Y)-\alpha_{14}(Y,X)=0.
\eeQ
Other terms vanishe by symmetry. Namely,
$$\alpha_1(X,Y)+\alpha_{10}(X,Y)+\alpha_{18}(X,Y)-\alpha_1(Y,X)-\alpha_{10}(Y,X)-\alpha_{18}(Y,X)=0,$$
$$\alpha_3(X,Y)+\alpha_{6}(X,Y)-\alpha_3(Y,X)-\alpha_{6}(Y,X)=0,$$
and
\beQ
\alpha_4(X,Y)+\alpha_{5}(X,Y)+\alpha_{15}(X,Y)+\alpha_{17}(X,Y)+\alpha_{19}(X,Y)\\
-\alpha_4(Y,X)-\alpha_{5}(Y,X)-\alpha_{15}(Y,X)-\alpha_{17}(X,Y)-\alpha_{19}(X,Y)=0.
\eeQ
The terms $\alpha_{2}$, $\alpha_{7}$, $\alpha_{8}$ et  $\alpha_{12}$ can be combined. Indeed, if we set $$\alpha=\alpha_{2}+\alpha_{7}+\alpha_{8}+\alpha_{12},$$ then
$$\begin{array}{lll}
\alpha(X,Y)-\alpha(Y,X)
&=&\eta^2\Big[ f^2\left( Y\cdot X-X\cdot Y\right) +Y\cdot T\cdot X\cdot T -X\cdot T\cdot Y\cdot T\Big]\cdot\varphi\\\\
&=&\eta^2\Big[ f^2\left( Y\cdot X-X\cdot Y\right) +||T||^2\left( Y\cdot X-X\cdot Y\right)\Big]\cdot\varphi\\\\
&& -2\eta^2\left( \left\langle X,T\right\rangle Y\cdot T-\left\langle Y,T\right\rangle X\cdot T\right) \cdot\varphi,
\end{array}$$
by taking $X$ and $Y$ as an orthonormal frame $\{e_1,e_2\}$, we have
$$\begin{array}{lll}
\alpha(e_1,e_2)-\alpha(e_2,e_1)
&=&2\eta^2\Big[-T_1e_2\left( T_1e_1+T_2e_2\right) +T_2e_2\left( T_1e_1+T_2e_2\right)\Big]\cdot\varphi\\
&&-2\eta^2e_1\cdot e_2\cdot\varphi\\\\
&=&-2\eta^2\omega\cdot\varphi+2\eta^2\left( T_1^2\omega\cdot+T_1T_2-T_1T_2+T_2^2\omega\cdot\right) \varphi\\\\
&=&-2\eta^2f^2\omega\cdot\varphi.
\end{array}$$
Always for $X=e_1$ and $Y=e_2$, 
$$\alpha_9(e_1,e_2)+\alpha_{13}(e_1,e_2)-\alpha_9(e_2,e_1)-\alpha_{13}(e_2,e_1)=2T\cdot\omega\cdot\varphi=-2J(T)\cdot\varphi,$$
where $J$ is the rotation of positive angle $\frac{\pi}{2}$ on $TM$.
Finally, we get
\beQ
\mathcal{R}(e_1,e_2)\varphi&=&\frac{1}{4}\left( Ae_2\cdot Ae_1-Ae_1\cdot Ae_2\right)\cdot\varphi-2\eta^2f^2\omega\cdot\varphi\\
&&-\frac{1}{2}d^{\nabla}A(e_1,e_2)\cdot\varphi-2J(T)\cdot\varphi. 
\eeQ
Using the Ricci identity 
$$\mathcal{R}(e_1,e_2)\varphi=-\frac{1}{2}R_{1212}e_1\cdot e_2\cdot\varphi,$$
we have
$$\left( \underbrace{R_{1212}-\ddet(A)-\kappa f^2}_{G}\right)\omega\cdot\varphi=\left( \underbrace{d^{\nabla}A(e_1,e_2)+\kappa^2fJ(T)}_{C}\right)\cdot\varphi,$$
that is,
$$G\omega\cdot\varphi=C\cdot\varphi,$$  
where $G$ is a function and $C$ a vector field. Since, $\omega\cdot\varphi=-i\overline{\varphi}$, we obtain
$$C\cdot\varphi^{\pm}=\pm iG\varphi^{\mp}.$$
Thus,
$$||C||^2\varphi^{\pm}=-G^2\varphi^{\pm}.$$
Since $\varphi$ is a non-trivial solution of (\ref{killinggeneral}), by Proposition \ref{normeconstante1}, $\varphi^+$ et $\varphi^-$ cannot vanish at the same point. Then $C=0$ and $G=0$. But $G=0$ is the Gauss equation and $C=0$ is the Codazzi equation. $\finpreuve$
\subsection{Generalized Killing Spinor and Dirac Equation}
\label{sect33}

If a spinor field $\varphi$ is a solution of the equation of generalized Killing spinors (\ref{killinggeneral}), then, we deduce immediately that $\varphi$ is also a solution of the corresponding Dirac equation
\beqt\label{diracgeneral}
D\varphi=(H-2\eta f)\varphi-\eta T\cdot\varphi.
\eeqt
In this section, we will see that this Dirac equation, which is weaker that the equation of generalized Killing spinors, is in fact equivalent for spinors with a norm satisfying the condition of Proposition \ref{normeconstante1}. For this, we begin by the following lemma
\begin{lem}\label{lienAQphi}
Let $\varphi$ be a non-trivial solution of (\ref{killinggeneral}). Then, the symmetric endomorphism $A$ is given by
$$A(X,Y)=2Q_{\varphi}(X,Y)+B(X,Y),$$ where $\qphi$ is the energy-momentum tensor associated with the spinor field $\varphi$, and $B$ is the symmetric tensor defined outside of the set of the zeros of $\varphi$ by
\beQ
B(X,Y)&=&-2\pre\left\langle \eta\left\langle X,Y\right\rangle T\cdot\varphi,\varphi\right\rangle -2\pre\left\langle \eta f\left\langle X,Y\right\rangle\varphi,\varphi\right\rangle\\
&&-\pre\left\langle  \eta\left( \left\langle X,T\right\rangle Y+ \left\langle Y,T\right\rangle X\right)\cdot\varphi,\varphi \right\rangle 
\eeQ
\end{lem}
We recall that the energy-momentum tensor associated with the spinor field $\varphi$ is the symmetric $2$-tensor defined by
$$\qphi(X,Y):=\frac{1}{2}\pre\left\langle X\cdot\nabla_Y\varphi+Y\cdot\nabla_X\varphi,\varphi/|\varphi|^2\right\rangle .$$
{\it Proof:}
A simple calculation of the energy-momentum tensor $\qphi$ using the fact that $\varphi$ is a solution of (\ref{killinggeneral}) yields to the result.
$\finpreuve$
\\
We can compute $B$ explicitely in both case, $\eta=\frac{1}{2}$ and $\eta=\frac{i}{2}$.
\begin{lem}\label{calculB}
\begin{itemize}
\item If $\eta=\frac{1}{2}$, then $B=-f\iid.$
\item If $\eta=\frac{i}{2}$, then $B$ satisfies
$$\left\lbrace 
\begin{array}{l}
B_{11}=-\pre\left\langle  iT\cdot\varphi,\varphi\right\rangle -\pre\left\langle iT_1e_1\cdot\varphi,\varphi\right\rangle \\ 
B_{12}=B_{21}=\frac{1}{2}\pre\left\langle  i(T_1e_2+T_2e_1)\cdot\varphi,\varphi\right\rangle \\ 
B_{22}=-\pre\left\langle  iT\cdot\varphi,\varphi\right\rangle -\pre\left\langle iT_2e_2\cdot\varphi,\varphi\right\rangle 
\end{array}\right. $$
\end{itemize}
\end{lem}
\noindent
{\it Proof:} The proof is immediate from Lemma \ref{lienAQphi}.$\finpreuve$\\
\indent
We recall that any spinor $\varphi$ can be decomposed in two parts under the action of the volume element $\omega_2$. Namely, $\varphi=\varphi^++\varphi^-$, with $\omega_2\cdot\varphi^{\pm}=\pm\varphi^{\pm}$. Then, we have
\beqt\label{Dphipm}
D\varphi^{\pm}=H\varphi^{\mp}-\eta T\cdot\varphi^{\pm}-2\eta f\varphi^{\mp}.
\eeqt
Now, we define the following endomorphisms
$$\qphi^{\pm}(X,Y)=\pre\left\langle \nabla_X\varphi^{\pm},Y\cdot\varphi^{\mp}\right\rangle,$$
and
\beQ
B^{\pm}(X,Y)&=&-\pre\left\langle \eta X\cdot T\cdot\varphi^{\pm},Y\cdot\varphi^{\mp}\right\rangle -\pre\left\langle \eta\left\langle X,T\right\rangle\varphi^{\pm},Y\cdot\varphi^{\mp} \right\rangle\\
 &&-\pre\left\langle \eta fX\cdot\varphi^{\mp},Y\cdot\varphi^{\mp}\right\rangle 
\eeQ 
On the other hand, we set
$$A^{\pm}=\qphi^{\pm}+B^{\pm},$$
and 
$$W=\frac{A^+}{|\varphi^-|^2}-\frac{A^-}{|\varphi^+|^2}.$$
From now, we will consider two cases, $\eta=\frac{1}{2}$ and $\eta=\frac{i}{2}$. We give all the details for the case $\eta=\frac{1}{2}$ and only specify minor differences for the case $\eta=\frac{i}{2}$.
\subsubsection{The case $\eta=\frac{1}{2}$}
\label{sect331}
In this section, we assume that the norm of $\varphi$ is constant and we will show that $W$ is identically zero. For that, we will show that $W$ is symmetric, trace-free with rank less or equal to 1. First, let us consider the trace.
\begin{lem}\label{traceQphi}
The endomorphisms $\qphi^{\pm}$ and $B^{\pm}$ satisfy
\begin{enumerate}
\item $\trace(\qphi^{\pm})=-H|\varphi^{\mp}|^2+\frac{1}{2}\pre\left\langle T\cdot\varphi^{\pm},\varphi^{\mp}\right\rangle +f|\varphi^{\mp}|^2.$
\item $\trace(B^{\pm})=-\frac{1}{2}\pre\left\langle T\cdot\varphi^{\pm},\varphi^{\mp}\right\rangle -f|\varphi^{\mp}|^2.$
\end{enumerate}
\end{lem}
\noindent
{\it Proof:} For $\qphi^{\pm}$, we have
\beQ
\trace(\qphi^{\pm})&=&\qphi^{\pm}(e_1,e_1)+\qphi^{\pm}(e_2,e_2)\\
&=&\pre\left\langle \nabla_{e_1}\varphi^{\pm},e_1\cdot\varphi^{\mp}\right\rangle +\left\langle \nabla_{e_2}\varphi^{\pm},e_2\cdot\varphi^{\mp}\right\rangle\\
&=&-\pre\left\langle D\varphi^{\pm},\varphi^{\mp}\right\rangle \\
&=&-H|\varphi^{\mp}|^2+\frac{1}{2}\pre\left\langle T\cdot\varphi^{\pm},\varphi^{\mp}\right\rangle +\pre\left\langle f\varphi^{\mp},\varphi^{\mp}\right\rangle.
\eeQ
For $B^{\pm}$, we have
\beQ
\trace(B^{\pm})&=&B^{\pm}(e_1,e_1)+B^{\pm}(e_2,e_2)\\
&=&-\pre\left\langle T\cdot\varphi^{\pm},\varphi^{\mp}\right\rangle+\frac{1}{2}\pre\left\langle T\cdot\varphi^{\pm},\varphi^{\mp}\right\rangle-\pre\left\langle f\varphi^{\mp},\varphi^{\mp}\right\rangle \\
&=&-\frac{1}{2}\pre\left\langle T\cdot\varphi^{\pm},\varphi^{\mp}\right\rangle -f|\varphi^{\mp}|^2.
\eeQ
$\finpreuve$
\indent
Now, we give the following lemma for the defect of symmetry of $\qphi^{\pm}$ and $B^{\pm}$.
\begin{lem}\label{defautsym}
The endomorphisms $\qphi^{\pm}$ and $B^{\pm}$ satisfy
\begin{enumerate}
\item $\qphi^{\pm}(e_1,e_2)=\qphi^{\pm}(e_2,e_1)+\frac{1}{2}\pre\left\langle \omega\cdot T\cdot\varphi^{\pm},\varphi^{\mp}\right\rangle ,$
\item $B^{\pm}(e_1,e_2)=B^{\pm}(e_2,e_1)-\frac{1}{2}\pre\left\langle \omega\cdot T\cdot\varphi^{\pm},\varphi^{\mp}\right\rangle.$
\end{enumerate}
\end{lem}
\noindent
{\it Proof:} We compute $\qphi^{\pm}(e_1,e_2)$:
\beQ
\qphi^{\pm}(e_1,e_2)&=&\pre\left\langle \nabla_{e_1}\varphi^{\pm},e_2\cdot\varphi^{\mp}\right\rangle \\
&=&\pre\left\langle e_1\cdot\nabla_{e_1}\varphi^{\pm},e_1\cdot e_2\cdot\varphi^{\mp}\right\rangle \\
&=&\pre\left\langle D\varphi^{\pm},\omega\cdot\varphi^{\mp}\right\rangle -\pre\left\langle e_2\cdot\nabla_{e_2}\varphi^{\pm},e_1\cdot e_2\cdot\varphi^{\mp}\right\rangle \\
&=&\pre\left\langle D\varphi^{\pm},\omega\cdot\varphi^{\mp}\right\rangle +\qphi^{\pm}(e_2,e_1).
\eeQ
From (\ref{Dphipm}), we get
\beQ
\pre\left\langle D\varphi^{\pm},\omega\cdot\varphi^{\mp}\right\rangle&=&\pre\left\langle H\varphi^{\mp},\omega\cdot\varphi^{\mp}\right\rangle-\frac{1}{2}\pre\left\langle  T\cdot\varphi^{\pm},\omega\cdot\varphi^{\mp}\right\rangle-f\pre\left\langle  \varphi^{\mp},\omega\cdot\varphi^{\mp}\right\rangle \\
&=&-\frac{1}{2}\pre\left\langle  T\cdot\varphi^{\pm},\omega\cdot\varphi^{\mp}\right\rangle\\
&=&\frac{1}{2}\pre\left\langle \omega\cdot T\cdot\varphi^{\pm},\varphi^{\mp}\right\rangle.
\eeQ
Now, we compute $B^{\pm}(e_1,e_2)-B^{\pm}(e_2,e_1)$.
\beQ
B^{\pm}(e_1,e_2)-B^{\pm}(e_2,e_1)&=&-\frac{1}{2}\pre\left\langle e_1\cdot T\cdot\varphi^{\pm},e_2\cdot\varphi^{\mp}\right\rangle -\frac{1}{2}\pre\left\langle T_1\varphi^{\pm},e_2\cdot\varphi^{\mp} \right\rangle\\
&&-\frac{1}{2}\pre\left\langle fe_1\cdot\varphi^{\mp},e_2\cdot\varphi^{\mp} \right\rangle+\frac{1}{2}\pre\left\langle e_2\cdot T\cdot\varphi^{\pm},e_1\cdot\varphi^{\mp}\right\rangle \\
&&+\frac{1}{2}\pre\left\langle T_2\varphi^{\pm},e_1\cdot\varphi^{\mp} \right\rangle+\frac{1}{2}\pre\left\langle fe_2\cdot\varphi^{\mp},e_1\cdot\varphi^{\mp} \right\rangle \\\\
&=&\frac{1}{2}\pre\left\langle e_2\cdot e_1\cdot T\cdot\varphi^{\pm},\varphi^{\mp}\right\rangle-\frac{1}{2}\pre\left\langle e_1\cdot e_2\cdot T\cdot\varphi^{\pm},\varphi^{\mp}\right\rangle\\
&&+\frac{1}{2}\pre\left\langle (T_1e_2-T_2e_1)\cdot\varphi^{\pm},\varphi^{\mp}\right\rangle\\\\
&=&- \pre\left\langle \omega\cdot T\cdot\varphi^{\pm},\varphi^{\mp}\right\rangle+\frac{1}{2}\pre\left\langle \omega\cdot T\cdot\varphi^{\pm},\varphi^{\mp}\right\rangle\\
&=&-\frac{1}{2}\pre\left\langle \omega\cdot T\cdot\varphi^{\pm},\varphi^{\mp}\right\rangle.
\eeQ
$\finpreuve$
From the last two lemmas, we deduce the following one.
\begin{lem}\label{symtrace}
The endomorphism $W$ is symmetric and trace-free.
\end{lem}
\noindent
{\it Proof:} 
From Lemma \ref{traceQphi} we deduce that
\beQ
\trace(A^{\pm})&=&\trace(\qphi^{\pm})+\trace(B^{\pm})\\
&=&-H|\varphi^{\mp}|^2+\frac{1}{2}\pre\left\langle T\cdot\varphi^{\pm},\varphi^{\mp}\right\rangle +f|\varphi^{\mp}|^2-\frac{1}{2}\pre\left\langle T\cdot\varphi^{\pm},\varphi^{\mp}\right\rangle -f|\varphi^{\mp}|^2\\
&=&-H|\varphi^{\mp}|^2
\eeQ
Thus
$$\trace(W)=\frac{\trace(A^+)}{|\varphi^-|^2}-\frac{\trace(A^-)}{|\varphi^+|^2}=0.$$
On the other hand, from Lemma \ref{defautsym}, we have
\beQ
A^{\pm}(e_1,e_2)&=&\qphi^{\pm}(e_1,e_2)+B^{\pm}(e_1,e_2)\\
&=&\qphi^{\pm}(e_2,e_1)+\frac{1}{2}\pre\left\langle \omega\cdot T\cdot\varphi^{\pm},\varphi^{\mp}\right\rangle +B^{\pm}(e_2,e_1)-\frac{1}{2}\pre\left\langle \omega\cdot T\cdot\varphi^{\pm},\varphi^{\mp}\right\rangle\\
&=&\qphi^{\pm}(e_2,e_1)+B^{\pm}(e_2,e_1)=A^{\pm}(e_2,e_1).
\eeQ
Then, $A^+$ and $A^-$ are symmetric and $W$ too.
$\finpreuve$\\
\noindent
Finally, we give this last lemma which show that $W$ has a rank at most 1.
\begin{lem}\label{rang}
The endomorphism field $W$ satisfies
$$\pre\left\langle W(X)\cdot\varphi^{-},\varphi^{+}\right\rangle =0$$
for any tangent vecor field $X$.
\end{lem}
\noindent
{\it Proof:}
First, since $\varphi$ has constant norm, for any tangent vector field $X$, we have
\beq\label{U}
0&=&X|\varphi|^2\nonumber\\
&=&X\left( |\varphi^+|^2+|\varphi^-|^2\right) \nonumber\\
&=&2\pre\left( \nabla_X\varphi^+,\varphi^+\right)+ 2\pre\left( \nabla_X\varphi^-,\varphi^-\right)\nonumber\\
&=&2\pre\left( U(X)\cdot\varphi^-,\varphi^+\right) ,
\eeq
where $U(X):=\dfrac{\qphi^+(X)}{|\varphi^-|^2}-\dfrac{\qphi^-(X)}{|\varphi^+|^2}$. In order to simplify the notations, we set
$$U^+(X):=\dfrac{\qphi^+(X)}{|\varphi^-|^2}\quad\text{and}\quad U^-(X):=\dfrac{\qphi^-(X)}{|\varphi^+|^2}.$$
On the other hand, we denote
$$V^+(X):=\dfrac{B^+(X)}{|\varphi^-|^2}\quad,\quad V^-(X):=\dfrac{B^-(X)}{|\varphi^+|^2}\quad\text{and}\quad V=V^+-V^-.$$
Let us compute $\pre\left( V(X)\cdot\varphi^-,\varphi^+\right)$.
\beQ
\pre\left( V(X)\cdot\varphi^-,\varphi^+\right)&=&\pre\left( V^+(X)\cdot\varphi^-,\varphi^+\right)-\pre\left( V^-(X)\cdot\varphi^-,\varphi^+\right)\\
&=&\pre\left( V^+(X)\cdot\varphi^-,\varphi^+\right)+\pre\left( V^-(X)\cdot\varphi^+,\varphi^-\right).
\eeQ
But ,
$$\pre\left( V^+(X)\cdot\varphi^-,\varphi^+\right)=\pre\left( V^+(X,e_1)e_1\cdot\varphi^-,\varphi^+\right)+\pre\left( V^+(X,e_2)e_2\cdot\varphi^-,\varphi^+\right).$$
Now, we compute the term $V^+(X,e_1)$.
\beQ V^+(X,e_1)&=&
-\frac{1}{2}\pre\left\langle X\cdot T\cdot\varphi^{+},e_1\cdot\frac{\varphi^{-}}{|\varphi^{-}|^2}\right\rangle-\frac{1}{2}\pre\left\langle X,T\right\rangle \left\langle \varphi^{+},e_1\cdot\frac{\varphi^{-}}{|\varphi^{\mp}|^2}\right\rangle \\
&&-\frac{1}{2}\pre\left\langle fX\cdot\varphi^{-},e_1\cdot\frac{\varphi^{-}}{|\varphi^{-}|^2}\right\rangle
\eeQ
Similarly, we have
 \beQ V^+(X,e_2)&=&
-\frac{1}{2}\pre\left\langle X\cdot T\cdot\varphi^{+},e_2\cdot\frac{\varphi^{-}}{|\varphi^{-}|^2}\right\rangle-\frac{1}{2}\pre\left\langle\left\langle X,T\right\rangle  \varphi^{+},e_2\cdot\frac{\varphi^{-}}{|\varphi^{-}|^2}\right\rangle \\
&&-\frac{1}{2}\pre\left\langle fX\cdot\varphi^{-},e_2\cdot\frac{\varphi^{-}}{|\varphi^{-}|^2}\right\rangle
\eeQ
Since $\left\{ e_1\cdot\frac{\varphi^-}{|\varphi^-|},e_2\cdot\frac{\varphi^-}{|\varphi^-|}\right\} $ is a local orthonormal frame of $\Gamma(\Sigma^+N)$ for the scalar product $\pre\left\langle .,.\right\rangle $, we deduce the following
$$V^+(X)\cdot\varphi^-=-\frac{1}{2}X\cdot T\cdot\varphi^{+}-\frac{1}{2}\left\langle X,T\right\rangle  \varphi^{+}-\frac{1}{2}fX\cdot\varphi^{-}.$$
Then, we conclude that 
\beQ
\pre\left\langle V^+(X)\cdot\varphi^-,\varphi^+\right\rangle &=& -\frac{1}{2}\pre\left\langle X\cdot T\cdot\varphi^+,\varphi^+\right\rangle-\frac{1}{2}\pre\left\langle \left\langle X,T \right\rangle \varphi^+, \varphi^+\right\rangle\\
&&-\frac{1}{2}\pre\left\langle fX\cdot\varphi^+,\varphi^- \right\rangle.
\eeQ
We can easily see that for any vector field $X$,
$$\pre\left\langle X\cdot T\cdot\varphi^+,\varphi^+\right\rangle+\pre\left\langle \left\langle X,T \right\rangle \varphi^+, \varphi^+\right\rangle=0,$$ which yields to
$$\pre\left\langle V^+(X)\cdot\varphi^-,\varphi^+\right\rangle=-\frac{1}{2}\pre\left\langle fX\cdot\varphi^+,\varphi^- \right\rangle.$$
By the same way, we show
$$\pre\left\langle V^-(X)\cdot\varphi^+,\varphi^-\right\rangle=-\frac{1}{2}\pre\left\langle fX\cdot\varphi^-,\varphi^+ \right\rangle.$$
Finally, we conclude that 
\beqt\label{V}
\pre\left\langle V(X)\cdot\varphi^{-},\varphi^{+}\right\rangle =0.
\eeqt
Since $W=U+V$, we deduce from (\ref{U}) and (\ref{V}) that 
$$\pre\left\langle W(X)\cdot\varphi^{-},\varphi^{+}\right\rangle =0.$$
This achieves the proof of the lemma.$\finpreuve$\\
The fact that
$$\pre\left\langle W(X)\cdot\varphi^{-},\varphi^{+}\right\rangle =0$$
implies that $W$ has rank at most 1. Since $W$ is also symmetric and trace-free, we deduce that 
\begin{prop}\label{W}
The endomorphisms field $W$ is identically zero.
\end{prop}
\noindent
Now, we can state the following proposition
\begin{prop}\label{diracimpliquekilling}
If $\varphi$ is a non-trivial solution of the Dirac equation (\ref{diracgeneral}) with constant norm, then $\varphi$ is also a solution of the equation of generalized Killing spinors (\ref{killinggeneral}).
\end{prop}
\noindent
{\it Proof:} We set $F:=A^++A^-$. The fact that $W=0$ implies 
\beqt\label{Fsurphi}
\frac{F}{|\varphi|^2}=\frac{A^+}{|\varphi^-|^2}=\frac{A^-}{|\varphi^+|^2}.
\eeqt
On the other hand, we have seen that for any vector field $X$, 
\beQ
\nabla_X\varphi&=&U^+(X)\cdot\varphi^-+U^-(X)\cdot\varphi^+\\
&=&\frac{A^+}{|\varphi^-|^2}\cdot\varphi^-+\frac{A^-}{|\varphi^+|^2}\cdot\varphi^+-V^+(X)\cdot\varphi^--V^-(X)\cdot\varphi^+,
\eeQ
which gives with (\ref{Fsurphi})
\beqt\label{nablaphi}
\nabla_X\varphi=\frac{F(X)}{|\varphi|^2}\cdot\varphi-V^+(X)\cdot\varphi^--V^-(X)\cdot\varphi^+.
\eeqt
The tensor $F$ is symmetric because $A^+$ and $A^-$ are symmetric. Then, it is easy to check that $A:=-2F$ is as in Lemma \ref{lienAQphi} and that (\ref{nablaphi}) gives
$$ \nabla_X\varphi=\frac{1}{2} X\cdot T\cdot\varphi+\frac{1}{2} fX\cdot\varphi+\frac{1}{2}\left\langle X,T\right\rangle \varphi-\frac{1}{2}AX\cdot\varphi.$$
This concludes the proof.
$\finpreuve$

\subsubsection{The case $\eta=\frac{i}{2}$}
\label{sect332}
The steps of the proof are exactly the same as for $\eta=\frac{1}{2}$. The two minor differences are that the expression of the spinor field solution of the Dirac equation is different and the norm of this spinor is not constant but satisfies
$$X|\varphi|^2= \pre\left\langle  iX\cdot T\cdot\varphi+ifX\cdot\varphi,\varphi\right\rangle. $$
We do not give the details but all the lemmas of last section have an analogue in the case $\eta=\frac{i}{2}$, and we obtain the following proposition.
\begin{prop}\label{diracimpliquekilling}
If $\varphi$ is a non-trivial solution of the Dirac equation (\ref{diracgeneral}) with constant norm, then $\varphi$ is also a solution of the equation of generalized Killing spinors (\ref{killinggeneral}).
\end{prop}
\noindent
Now, we have all ingredients to prove Theorems \ref{thm1} and \ref{thm2}.
\subsection{Proof of Theorems \ref{thm1} and \ref{thm2}}
The proof is the same for the two cases.\\
The equivalance between $(i)$ and $(ii)$ has been proved in Section \ref{sect32}. On the other hand, it is clear that $(ii)$ implies $(iii)$. Finally, we have to show that $(iii)$ implies $(ii)$. We have seen in Section  \ref{sect33} that a spinor field solution of the Dirac equation (\ref{diracgeneral}) with constant norm is also solution of the equation of generalized Killing spinors (\ref{killinggeneral}). We just have to show that
$$\nabla_XT=fAX$$
and
$$
df(X)=-\left\langle AX,T\right\rangle .$$
The second identity is obvious since $$df=-2Q_{\varphi}(T)-B(T)$$ and the expression of $A$.\\
 To show the first point, we set $$a(X):=\nabla_XT-fAX.$$
 Since $\langle\nabla_XT,Y\rangle=\langle\nabla_YT,X\rangle$, it is clear that $a$ is symmetric. Moreover, since $2Hf=\ddiv(T)$, we deduce that $a$ is trace-free. Finally, $\nabla_TT=-f\diff f$ implies the following relation
 $$\nabla_TT=-2Q_{\varphi}(T)-B(T),$$
 which is equivalent to $a(T)=0$. Consequently, $a$ has rank at most 1. So, we can conlude that $a$ is identically zero. It means that for any vector field $X$ on $N$, we have
 $$\nabla_XT=fAX.$$
 This achieves the proof.$\finpreuve$
 
\section{Isometric immersions into  $\Ekt$, $\tau\neq0$}
Here, we state the main results about spinorial characterization of surfaces isometrically immersed into $\Ekt$.
\begin{thm}\label{thm3} 
Let $(N,\langle.,.\rangle)$ be a connected, simply connected and oriented Riemannian surface. Let $T$ be a vector field, $f$ and $H$ two functions on $N$ satisying
$$\left\lbrace  \begin{array}{l}
f^2+||T||^2=1, \\ 
\langle\nabla_XT,Y\rangle=\langle\nabla_YT,X\rangle+2\tau f\left\langle X,JY\right\rangle ,\quad\forall X,Y\in\chi(N), \\ 
2Hf=\ddiv(T).
\end{array}\right. $$
Let $A\in{\rm{End}}(TN)$. The following three statements are equivalent.
\begin{enumerate}[i)]
\item There exists an isometric immersion $F$ from $N$ into $\Ekt$ with mean curvature $H$ such that the Weingarten operator related to the normal $\nu$ is given by
$$dF\circ A\circ dF^{-1}$$ 
and such that
$$\xi=dF(T)+f\nu.$$
\item The spinor field $\varphi$ is solution of the equation
$$\nabla_X\varphi=-\frac{\tau}{2} X\cdot\omega\cdot\varphi+\frac{\alpha}{2}\langle X,T\rangle T\cdot\omega\cdot\varphi-\frac{\alpha}{2}
f\langle X,T\rangle \omega\cdot\varphi-\frac{1}{2}AX\cdot\varphi,$$
where $A$ satisfies 
$$\nabla_XT=f(AX-\tau JX),$$
and
$$df(X)=-\langle AX-\tau JX,T\rangle.$$
\item The spinor field $\varphi$ is solution of the Dirac equation 
$$D\varphi=H\varphi+\left( \tau-\frac{\alpha}{2}||T||^2\right) \omega\cdot\varphi-\frac{\alpha}{2}
fT\cdot\omega\cdot\varphi,$$ with constant norm and such that
$$df=-2Q_{\varphi}(T)-B(T)-\tau J(T),$$
where $B$ is the tensor defined in an orthonormal frame $\{e_1,e_2\}$ by the following matrix
$$
\left( 
\begin{array}{cc}
\alpha T_1T_2 & \frac{\alpha}{2}\left( T_1^2-T_2^2\right) \\ \\
\frac{\alpha}{2}\left( T_1^2-T_2^2\right) & -\alpha T_1T_2
\end{array}
\right) ,
$$
with $T_i=\left\langle T,e_i\right\rangle $.
\end{enumerate}
\end{thm}
\pf
The proof of Theorem \ref{thm3} uses the same arguments in the proof of Theorems \ref{thm1} and \ref{thm2}, so we do not give it here. Just note that we use another theorem of Daniel (\cite{Dan2}) which is the analogue of the theorem we give in Section \ref{sect31} for isometric immersions inot $\Ekt$.
\begin{rem}
 In \cite{Dan2}, Daniel gives a characterization of surfaces isometrically immersed into 3-homogeneous manifolds with 4-dimensional isometry group which are not a product. Then, he deduces a Lawson correspondence for cmc-surfaces into such manifolds and product spaces $\M^2(\kappa)\times\R$. The Lawson correspondence is a local isometric correspondence between constant mean curvature surfaces into homogeneous 3-manifolds. Then, it would be natural to understand the Lawson correspondence in terms of spinors. Note that the classical Lawson correspondence (for cmc-surfaces in space forms) is still not well understood via spinors.
 \end{rem}
\noindent
\textbf{Acknowledgement.} The author would like to thank Oussama Hijazi and Jean-Fran\c{c}ois Grosjean for their encouragement and helpful remarks.

\bibliographystyle{amsplain}
\bibliography{mabiblio}
\end{document}